\documentclass[a4paper,10pt]{article}

\usepackage[english]{babel}
\usepackage{german, umlaute, amsmath, amsfonts,
amssymb , latexsym, theorem, epic, epsfig, rotating, graphicx,enumerate}

\theoremstyle{break}
\newtheorem{definition}{Definition}[section]

\newtheorem{lem}[definition]{Lemma}
\newtheorem{coro}[definition]{Corollary}
\newtheorem{cor}[definition]{Corollary}
\newtheorem{prop}[definition]{Proposition}

\newtheorem{lemma}[definition]{Lemma}
\newtheorem{Def}[definition]{Definition}

\theoremstyle{plain}
\newtheorem{rmk}[definition]{Remark}

\theorembodyfont{\rmfamily}

\numberwithin{equation}{section}

\newcommand{\de}{\delta}  
 
\newcommand{\e}{\varepsilon} 
\newcommand{\om}{\omega} 
 
\newcommand{\gt}{\theta} \newcommand{\BR}{\mathbb{R}}
\newcommand{\BT}{\mathbb{T}} 

\newcommand{\BZ}{\mathbb{Z}} \newcommand{\BN}{\mathbb{N}}
\newcommand{\A}{\mathcal{A}}

 \newcommand{\sdic}{{\em sdic}}

\newcommand{\OC}[1]{\overline{O(#1)}}
\newcommand{\SD}{\text{\sl SD}}
\newcommand{\LS}{\text{\sl LS}}
\newcommand{\Tr}{\text{\sl Tr}}
\newcommand{\diam}{\operatorname{diam}}
\newcommand{\card}{\operatorname{card}}

\newcommand{\J}{{\mathcal J}}
\newcommand{\rz}{{\mathbb R}}
\newcommand{\nz}{{\mathbb N}}
\newcommand{\tz}{{\mathbb T}}
\newcommand{\zz}{{\mathbb Z}}
\newcommand{\st}{{\rm\ s.t.\ }}
\newcommand{\eg}{{e.g.\ }}
\newcommand{\ie}{{i.e.\ }}
\newcommand{\qed}{\phantom{xxxxxx}\hfill q.e.d.\\}
   \newenvironment{proof}{{\em Proof:} }{\qed}

\newcommand{\xtil}{\ensuremath{\tilde{x}}}
\newcommand{\Lambdatil}{\ensuremath{\tilde{\Lambda}}}
\newcommand{\ytil}{\ensuremath{\tilde{y}}}
\newcommand{\ztil}{\ensuremath{\tilde{z}}}
\newcommand{\mutil}{\ensuremath{\tilde{\mu}}}

\newcommand{\twonorm}[1]{\ensuremath{\parallel #1 \parallel_{2}}}
\newcommand{\Matrix}[2]{\ensuremath{\left(\begin{array}{#1} #2
     \end{array}\right)}}

\title{\textsc{How chaotic are strange nonchaotic attractors?}}  \author{Paul
  Glendinning,\thanks{School of Mathematics, University of Manchester, Oxford
    Road, Manchester M13 9PL, United Kingdom. E-mail: {\tt
      p.a.glendinning@manchester.ac.uk}}\ \ and Tobias H. J\"ager and Gerhard
  Keller\thanks{both Mathematisches Institut, Universit\"at
    Erlangen-N\"urnberg, Bismarckstr.  $1\frac{1}{2}$, 91054 Erlangen,
    Germany. Email: jaeger@mi.uni-erlangen.de, keller@mi.uni-erlangen.de}}

\begin{document}
\maketitle

\begin{abstract}
  We show that the classic examples of quasi-periodically forced maps
  with strange nonchaotic attractors described by Grebogi {\em et al}
  and Herman in the mid-1980s have some chaotic properties. More
  precisely, we show that these systems exhibit sensitive dependence
  on initial conditions, both on the whole phase space and restricted
  to the attractor. The results also remain valid in more general
  classes of quasiperiodically forced systems. Further, we include an
  elementary proof of a classic result by Glasner and Weiss on
  sensitive dependence, and we clarify the structure of the attractor
  in an example with two-dimensional fibers also introduced by Grebogi
  {\em et al}.
\end{abstract}

\section{Introduction}

Strange nonchaotic attractors (SNA) are attractors of dynamical systems
which have some form of local contraction, but which also have a
complicated or fractal structure (hence the word `strange'). In the
context of quasiperiodically forced maps, i.e. maps of the form
\begin{equation}\label{qp}
f(\theta,\xi)=(\theta +\om ,f_\theta(\xi))
\end{equation}
for $\theta\in\BT$, $\om$ irrational, $\xi$ in some suitable metric
space $M$ and fiber maps $f_\theta$ defined by $f_\theta(\xi) =
\pi_2\circ f(\theta,\xi)$, this notion is used for compact invariant
sets which are the topological closure of a non-continuous invariant
graph with negative Lyapunov exponents in the fibers (see
Remark~\ref{rmk:sna} for the precise definition). The
negative Lyapunov exponent in the $\xi -$direction%
\footnote{The Lyapunov exponent in the $\theta-$direction associated
  to the corresponding invariant measure is always zero.}
provides local contraction and the topological entropy of the system
is zero.  These two conditions are generally considered sufficient to
justify calling the attractors nonchaotic, and the first one ensures
that there is local exponential convergence to the attractor in almost
all fibers of constant $\theta$.  Moreover, many authors remark that
this implies that there is no exponential sensitivity to initial
conditions. It has also been observed that the existence of a SNA
implies that finite time Lyapunov exponents can be positive \cite{PF},
and that this chaotic-like property is responsible for the lack of
smooth invariant curves.

In this paper we consider the chaotic-like properties of SNA in
quasiperiodically forced systems in more detail. We focus on the
property of sensitive dependence on initial conditions (\sdic ), which
has been regarded as one of the hallmarks of chaos, and show that many
SNA have this property.  By \sdic~ we mean the classic topological
definition as in Devaney \cite{D} which does not impose any conditions
on rates of separation. To be precise, the standard definition is as
follows.
\begin{Def}
Let $X$ be a metric space with metric $d$. A map $g:X\to X$ has
sensitive dependence on initial conditions (\sdic ) iff there exists
$\e >0$ such that for all $\de >0$ and $x\in X$ there exists $n\ge
0$ and $y\in X$ depending on $x$ and $\de$ such that
$d(x,y) <\de$ and $d(g^n(x),g^n(y)) > \e$.
\end{Def}
This was one of the three conditions for chaos introduced by Devaney
\cite{D}, although it was later shown that it is implied by the other
two conditions (transitive and dense periodic orbits) \cite{Banks,GW1}
and so, as Glasner and Weiss observe \cite{GW1}, Devaney's definition
is too weak to be considered as a good definition of chaotic dynamics.
On the other hand it is certainly a feature associated with chaos, and
the presence of this property in SNAs emphasizes their position on the
cusp between regular and chaotic systems. It is also worth noting that
quasiperiodically forced systems cannot be chaotic in the sense of
Devaney because since $\om$ is irrational there are no periodic
orbits.  The definition of \sdic\ given above leaves some latitude in
the choice of the space $X$. First of all, it is natural to consider
the dynamics restricted to the attractor, thus choosing $X=\A$, and
this is investigated in section three.  Another obvious choice is to
set $X$ to be the whole space on which the system is defined, i.e.\
$X=\BT \times M$, and this will be treated in section four. The
difference between these choices is reflected in changes in the set of
points in a neighbourhood of any point.

It is one of the most interesting aspects of SNAs that the
measure-theoretic and topological point of view often separate, and
properties which are generic in the one sense are degenerate in the
other and vice versa. For example, it usually makes a great difference
whether the measure-theoretic or the topological support of an
invariant measure is considered, and there are situations where the
former has a very complicated structure while the later is just a
smooth torus.  In order to fully understand the behavior of
quasiperiodically forced maps it is often necessary to combine both
viewpoints, and this is explains that while the focus of this paper
lies on the topological side, measures will inevitably make an
appearance. (On a technical level this happens via the results of Glasner
and Weiss \cite{GW1}, for which we included a version of the proof,
with a strongly simplified measure-theoretic part in section two.)

One of the most studied classes of SNA arise in pinched skew products.  These
are systems (\ref{qp}) for which there exists at least one value of $\gt$,
$\gt^*$ say, such that $f(\theta^*,\xi)=0$ for all $\xi\in M$.  In other
words, at least one fiber of constant $\theta$ is mapped to a single point,
the pinched point. These systems include some of the original examples
suggested by \cite{GOPY}, and are one of the few classes of systems for which
it is possible to prove rigorous results about the existence and structure of
SNAs \cite{K,G,J1}.  It is not hard to adapt the results of \cite{GW1} to
pinched skew products which satisfy three natural conditions to prove:
\begin{quote}
\textit{If $\A$ is the attractor of a pinched skew product
   $f:\BT\times \BR\to\BT\times \BR$ satisfying the conditions
   \eqref{property:1}\,--\,\eqref{property:3} of section three
   and $\A$ is not a continuous graph then $f$ has \sdic~ on $\A$. In
particular,
   if $\A$ is a SNA of a pinched skew product then $f$ has \sdic{} on $\A$.}
\end{quote}
See Corollary \ref{coro:1}. The results of Glasner and Weiss \cite{GW1}
give even more information
about the structure of points in an SNA of a pinched skew product. A
point $x\in X$ is {\em Lyapunov stable} for the map $g:X\to X$ if for
all $\e >0$ there exist $\de >0$ such that $d(g^n(x),g^n(y))<\e$ for
all $n\ge 0$ and $y\in X$ with $d(x,y)<\de$.  The existence
or nonexistence of Lyapunov stable points effectively determines the
dynamics of the pinched skew product.  \textit{
\begin{quote}
If $\A$ is the attractor of a pinched skew product which satisfies the
conditions \eqref{property:1}\,--\,\eqref{property:3} of section three then the
following are equivalent
\begin{enumerate}
\item $\A$ contains a Lyapunov stable point;
\item $\A$ is a continuous graph;
\item $\A$ does not have \sdic .
\end{enumerate}
\end{quote}
}
\noindent
This can be reinterpreted as saying that $\A$ contains no
Lyapunov stable points if and only if $\A$ is strange if and only if
$\A$ has \sdic . These statements are direct
consequences of more general results in section three: Corollary
\ref{coro:1} proves that (2) and (3) are equivalent, and the equivalence
of (1) then follows from the dichotomy (\ref{eq:dichotomy}). Note that the two
results stated above hold for the SNA in the classic example of Grebogi {\em et
al}
\cite{GOPY} which has $M=\BR$ and
\begin{equation}
f_\theta (\xi )= B\cos 2\pi\theta \tanh \xi
\end{equation}
in (\ref{qp}). This has a pinched SNA if $B>2$ \cite{GOPY,K}.

In section four we turn to the question of \sdic~on the whole phase space, using
techniques based only on the dynamics of quasiperiodically forced
one-dimensional maps. For the case of pinched skew products, we thus obtain
\sdic\
on the whole phase space whenever the attractor is not a continuous
graph, similar to the results above.

It appears harder to prove the existence of SNA in non-pinched cases.
The most prominent and for a long time also the only class of
quasiperiodically forced systems where this was possible are
quasiperiodic matrix cocycles \cite{herman:1983}, with quasiperiodic
Schr\"odinger cocycles as a special case. Only recently more general
approaches have been developed which are at least in principle
applicable to a much broader class of systems
(\cite{bjerkloev:2003,bjerkloev:2005,jaeger:2005b}) and thus confirm
the strong numerical evidence for the widespread existence of SNA in
quasiperiodically forced maps.  The application of our results to
these examples is discussed in more detail in sections three and four.
One important concept in this context is the rotation number of a
quasiperiodically forced circle homeomorphism, which Herman has shown
to exist in \cite{herman:1983} and which has been investigated further by
many authors. In section 4 we discuss convergence properties of the
rotation number and their implications for \sdic.

\begin{figure}[h]\label{fig:1}
      \begin{center}
        \includegraphics[width=7cm]{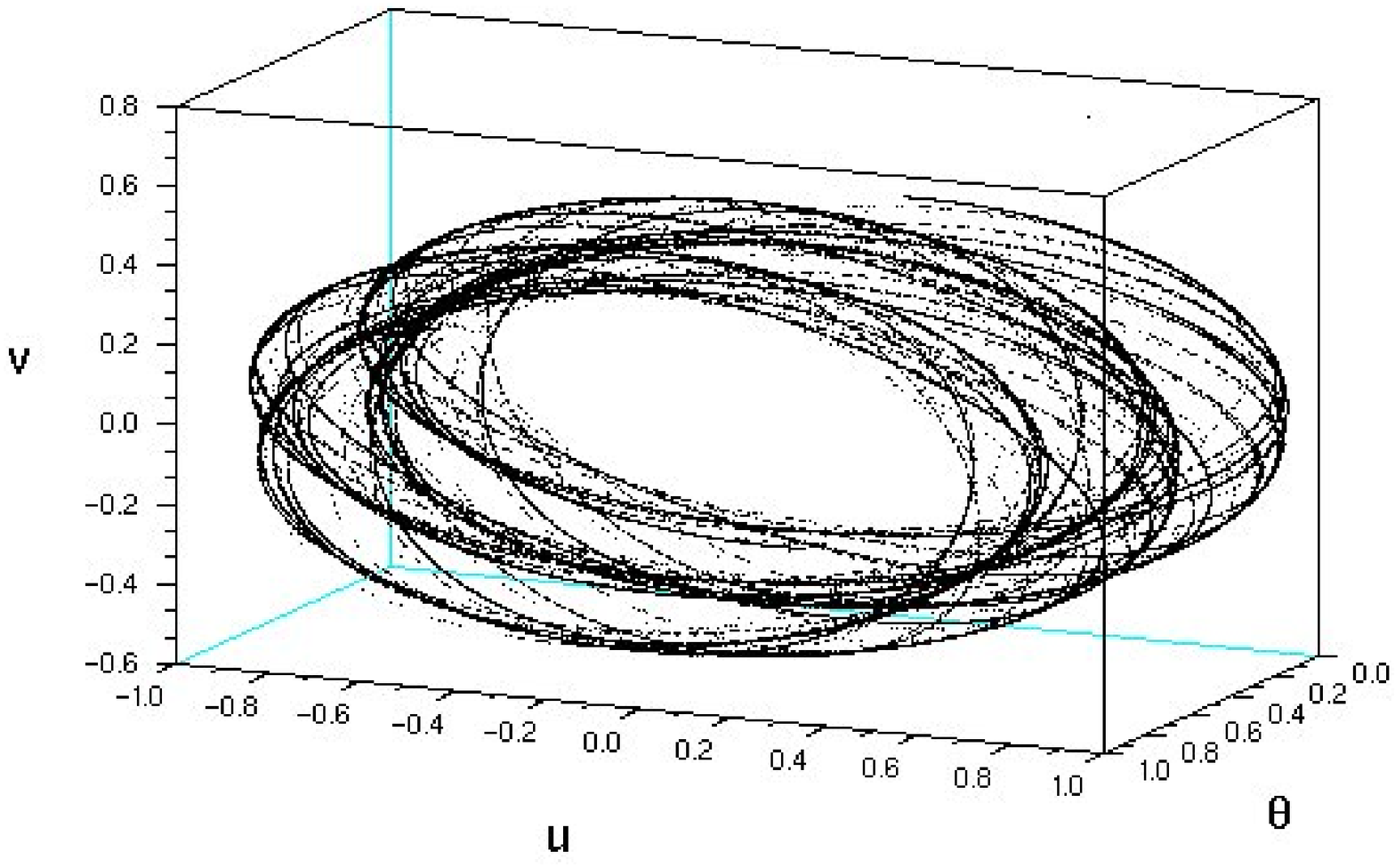}
        \hspace*{0.5cm}
        \includegraphics[width=7cm]{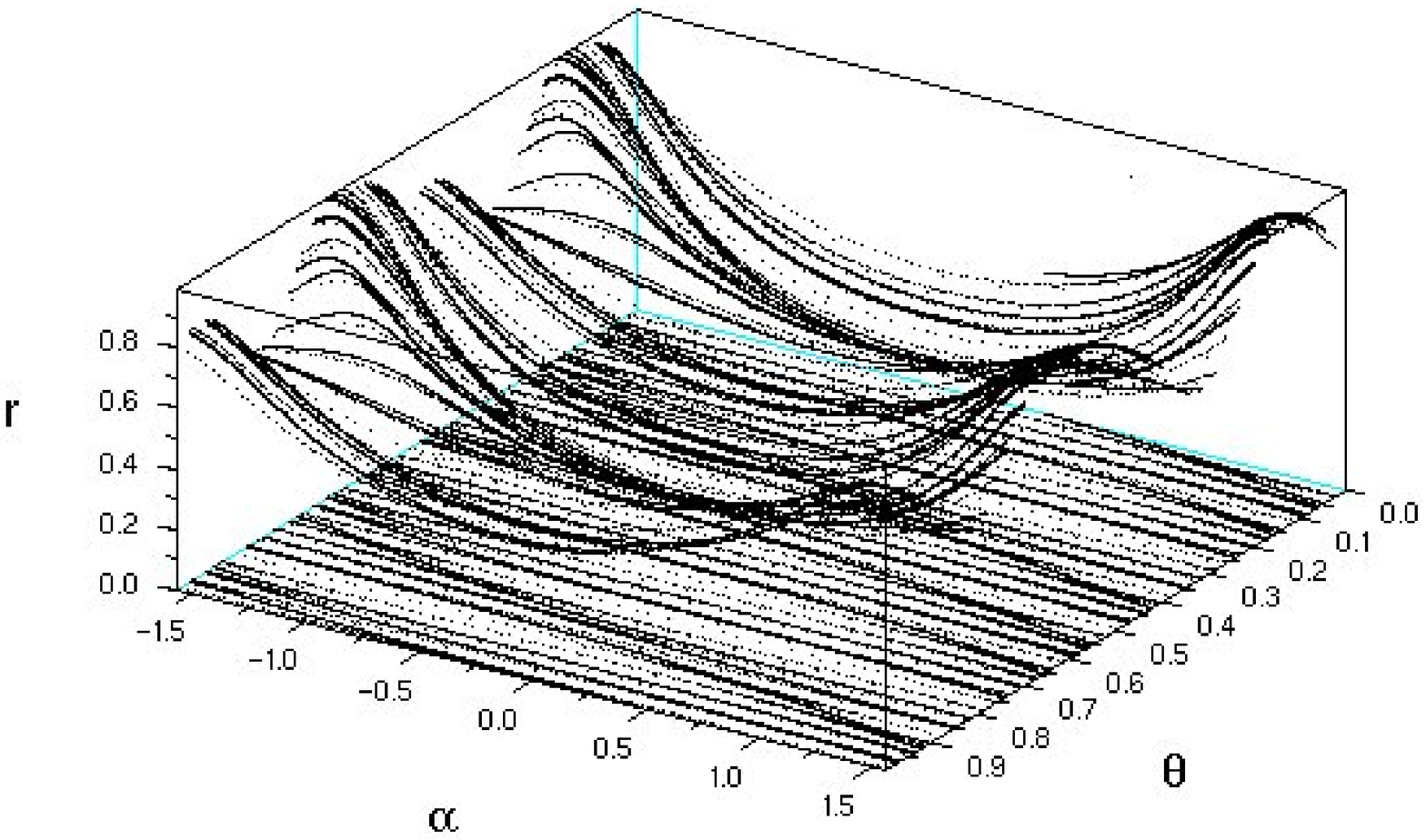}
        \caption{The figure to the left shows the attractor of the map
          $\Lambda:\tz1\times\rz2\to\tz1\times\rz2$, $
          \Lambda(\theta,u,v) := \big(\theta + \omega,
          \frac{\beta}{1+u2+v2} \,{1\ 0\choose 0\ \gamma} \,
          R_\theta\, {u\choose v}\big)$ where $R_\theta$ is the
          rotation matrix with angle $2\pi\theta$.  The figure to the
          right shows the same attractor when projective polar
          coordinates $(\alpha,r)$ are used in the fibers. It is
          plotted together with its projection to the
          $(\theta,\alpha)$-plane, which arises as the attractor of a
          quasiperiodically matrix cocycle of the type discussed by
          Herman. Details are given in Section~\ref{sec:final}.  }
      \end{center}
    \end{figure}

Finally, in section 5 we return to the second example from the
original paper of Grebogi et al \cite{GOPY}, which does not appear to
have been considered further in the literature so far. This system has
two-dimensional fibers and a non-pinched attractor, such that our
previous results do not apply directly. But after passing to
projective polar coordinates in the fibers we are able to clarify the
structure of the attractor and to relate its dynamics to those of a
matrix cocycle, which makes it possible to prove \sdic\ both on the
attractor and on the whole phase space. See Figure~\ref{fig:1}.

The existence of \sdic~ in SNAs will not come as a complete surprise
to the more applied community, although there has clearly been some
confusion.  Pikovsky and Feudel \cite{PF} define a quantity which
measures local separation due to changes in $\theta$ for
a given orbit which is a function
of the orbit and the number of iterates, $N$, on which separation is
considered.

Since SNA have non-positive Lyapunov exponents
their measure of separation cannot increase exponentially, but their careful
numerical
experiments suggest that the maximum separation over $N$
iterates grows roughly linearly with $N$ (to be more precise, their
experiments give a growth rate of $N^\mu$ with $\mu\approx 0.97$
\cite{PF}). In some sense our results can be seen as confirming that
their phase sensitivity exponent reflects \sdic~ in the system. Of course, for
forced differential equations this implies sensitive dependence with
respect to
small changes in the initial time of a solution as well as with respect to the
phase space.

\textit{Notation:} We reserve the letter $f$ to denote quasiperiodic maps of
the form (\ref{qp}), often with the choice of $M$ fixed to be $\rz$ or $\tz$.
General maps of a metric space $X$ will usually be denoted by $g$ (as in this
introduction). In particular, the results of section two are all in terms of
such general maps $g$.

\section{Sensitivity and equicontinuity}
\label{sec:dichotomy}

Let $g:X\to X$ be a continuous map of a compact metric space $(X,d)$ which
has no isolated points.  For $x\in X$ and $r>0$ let $B_r(x):=\{y\in
X:d(x,y)\leq r\}$ and denote by $\OC{x}$ the closure of the orbit
$\{g^n(x):n\in\nz\}$ of $x$. The set of
\emph{transitive} points, \ie the set of points $x$ for which $\OC{x}=X$ is
denoted by $\Tr$. If $\Tr\neq\emptyset$ one says that $g$ is transitive.

We are interested in the \emph{sensitive} and in the \emph{Lyapunov stable}
points of $g$. To this end we introduce, for each $\epsilon>0$, the two sets
\begin{align}
\label{eq:def-SDe}
   \SD_\epsilon&:=\{x\in X: \forall\delta>0\,\exists y,z\in
B_\delta(x)\,\exists
   n\in\nz \st d(g^ny,g^nz)\geq \epsilon\}\\
   \label{eq:def-LSe}
   \LS_\epsilon&:=\{x\in X: \exists\delta=\delta(\epsilon,x)>0\,\forall y,z\in
B_\delta(x)\,\forall
   n\in\nz: d(g^ny,g^nz)<\epsilon\}
\end{align}
Evidently $\LS_\epsilon=X\setminus\SD_\epsilon$. Let
\begin{displaymath}
\LS:=\bigcap_{\epsilon>0}\LS_\epsilon\text{ and }
\SD:=\bigcup_{\epsilon>0}\SD_\epsilon\ ,\text{ so }\LS=X\setminus\SD\;.
\end{displaymath}
$\LS$ is the set of \emph{Lyapunov stable} points, $\SD$ that of
\emph{sensitive} points. One says that the map $g$ has \emph{sensitive
dependence}, if $\SD_\epsilon=X$ for some $\epsilon>0$.\footnote{It
is easily seen that this definition is equivalent to the one given
in Def~1 in the introduction.} In that case, each
point of $X$ is sensitive, but the converse is not necessarily true.
However,
it follows immediately from these definitions that each $\SD_\epsilon$
is closed
and
forward invariant under $g$. Therefore, if $\SD$ contains a transitive
point $x$, then $X=\OC{x}\subseteq\SD_\epsilon$ for some $\epsilon>0$.

On the other hand, if no sensitive point is transitive, \ie if
$\Tr\subseteq\LS$, and if $\Tr\neq\emptyset$, then actually $\Tr=\LS$.%
\footnote{Let $y\in\Tr$, $x\in\LS_\epsilon$ for some $\epsilon>0$, and let
  $z\in X$ be arbitrary. There is $k\in\nz$ such that
  $d(g^ky,x)<\delta(\frac\epsilon2,x)$, $\delta$ as in \eqref{eq:def-LSe}.
  As $X=\OC y$ has no isolated point, there is $n>k$ such that
  $d(g^ny,z)<\frac\epsilon2$. Hence

  $d(g^{n-k}x,z)\leq d(g^{n-k}x,g^{n-k}(g^ky))+d(g^ny,z)\leq
\frac\epsilon2+\frac\epsilon2 =\epsilon$.  Therefore $z\in\OC x$. }
Hence we have the following dichotomy for transitive systems:%
\footnote{ Note that in transitive systems the set $\Tr$ is \emph{residual},
  \ie it contains a countable intersection of dense open sets. In particular
  it is dense by Baire's category theorem.}
\begin{equation}
\label{eq:dichotomy}
\text{If $g$ is transitive, then either $g$ has sensitive dependence or
$\Tr=\LS$. }
\end{equation}
Note that if $\LS=X$, then the family
$(g^n)_{n\in\nz}$ is actually \emph{equicontinuous}, \ie
\begin{equation}\label{eq:equicontinuous-def}
\forall\epsilon>0\,\exists\delta>0\,\forall y,z\in X\,\forall n\in\nz:\
d(y,z)\leq\delta\Rightarrow d(g^ny,g^nz)\leq\epsilon\;.
\end{equation}
This is an immediate consequence of the compactness of $X$.  So, if
$g$ is minimal (\ie if $\Tr=X$), then either $g$ has sensitive
dependence or $(g^n)_{n\in\nz}$ is equicontinuous, see \eg
\cite{AAB-96} where these and related questions are treated
systematically.  Glasner and Weiss \cite{GW1} showed that this
dichotomy remains true if the assumption of minimality is replaced by
the weaker one that $f$ is transitive and admits a finite invariant
measure with full topological support.\footnote{They also provide
  examples showing that transitivity alone is not sufficient for the
  dichotomy.}  In the rest of this section we will rederive this
result with a completely elementary self-contained proof that does
neither use Birkhoff's ergodic theorem nor any knowledge about
syndetic sets as does the proof in \cite{GW1}.

Recall that a point $x$ is \emph{nonwandering}, if
\begin{displaymath}
\label{eq:nonwandering}
\forall\delta>0\,\exists n>0 \st B_\delta(x)\cap
g^n(B_\delta(x))\neq\emptyset\;.
\end{displaymath}
It follows at once that, if $x$ is nonwandering, then for each $\delta>0$ the
set
\begin{displaymath}
\label{eq:Rdeltax}
R_\delta(x):=\{n>0: B_\delta(x)\cap
g^n(B_\delta(x))\neq\emptyset\}\quad\text{is infinite.}
\end{displaymath}
Observe also that, by definition, each transitive point is
nonwandering. Recall further that $g :X \rightarrow X$ is called
\emph{uniformly rigid} if there exists a sequence $n_k$ of integers
going to infinity, such that $g^{n_k}$ converges uniformly to the
identity map on $X$. Obviously, every uniformly rigid map must be a
homeomorphism.
\begin{lemma} \label{lem:uniformlyrigid}
If $x$ is a nonwandering Lyapunov point with $\delta=\delta(\epsilon/2,x)$ as
in \eqref{eq:def-LSe}, then
\begin{displaymath}
   \forall\epsilon>0\,\forall n\in R_{\delta(\epsilon/2,x)}(x)\,\forall
y\in\OC{x}:\
   d(g^ny,y)\leq\epsilon\ .
\end{displaymath}
In particular, $g|_{\OC x}$ is \emph{uniformly
   rigid}.
\end{lemma}
\begin{proof}
Let $B=B_{\delta(\epsilon/2,x)}$ and $n\in R_{\delta(\epsilon/2,x)}(x)$. Let
$u\in
B\cap g^nB$, $v\in B\cap g^{-n}\{u\}$. Then $v,g^n(v)\in
B=B_{\delta(\epsilon/2,x)}$ so that, for all $k\in\nz$,
\begin{displaymath}
   d(g^{k+n}x,g^kx)
   \leq
  d(g^{k+n}x,g^{k+n}v)+d(g^{k}(g^nv),g^kx)
   \leq
   \frac\epsilon2+\frac\epsilon2
   =
   \epsilon\;.
\end{displaymath}
Hence $d(f^ny,y)\leq\epsilon$ for all $y\in\OC x$.
\end{proof}

\noindent
To get more out of this one needs to control the sets $R_\delta(x)$ in
the previous lemma. More precisely, one needs to make sure that
$R_\delta(x)$ is \emph{syndetic}, \ie\ has bounded gaps:
\begin{equation}
\label{eq:bounded-gaps}
\exists s\in\nz\,\forall n\in\nz:\ [n,n+s]\cap R_\delta(x)\neq\emptyset\;.
\end{equation}
In this case we say that the gaps are bounded by $s$.
\begin{lemma}\label{lem:g2}
If $x$ is a nonwandering Lyapunov point for which all sets $R_\delta(x)$
$(\delta>0)$ have bounded gaps, then the family $(g^n|_{\OC x})_{n\in\nz}$
is equicontinuous.
\end{lemma}
\begin{proof}
Let $\epsilon>0$ and choose $\delta=\delta(\epsilon/6,x)$ as in
\eqref{eq:def-LSe}. We may assume that $\delta\in(0,\epsilon/3)$. Let the
gaps of $R_\delta(x)$ be bounded by $s$. As all $f^j$ are continuous, there
is some $\eta\in(0,\delta]$ such that
\begin{displaymath}
   \forall y,z\in X\,\forall j\in\{0,\dots,s\}:\
d(y,z)\leq\eta\Rightarrow d(g^jy,g^jz)<\epsilon/3\;.
\end{displaymath}
Let $y,z\in\OC x$.  Any $n$ can be written as $n=k+j$ with $k\in R_\delta(x)$
and $j\in\{0,\dots,s\}$. Hence, if $d(y,z)\leq\eta$ then, by
Lemma~\ref{lem:uniformlyrigid},
\begin{displaymath}
   d(g^ny,g^nz)
   \leq
   d(g^k(g^jy),g^jy)+d(g^jy,g^jz)+d(g^jz,g^k(g^jz))
   <
   \epsilon\;.
\end{displaymath}
\end{proof}

\noindent
It remains to give a condition ensuring that the sets
$R_\delta(x)$ have bounded gaps. Since this is a kind of uniform recurrence
condition, the following lemma is not too surprising.
\begin{lemma}\label{lem:g3}
If $x$ is a nonwandering Lyapunov point which belongs to the
topological support of a $g$-invariant
finite measure $\mu$ on $X$\footnote{This
   means $\mu(U)>0$ for each open neighbourhood $U$ of $x$.}, then
$R_\delta(x)$ has bounded gaps for each $\delta>0$.
\end{lemma}
\begin{proof}
  As $\mu$ has full topological support, the ball $B:=B_{\delta}(x)$
  has positive $\mu$-measure for every $\delta>0$. Denote $\hat
  B:=\bigcup_{j=0}^\infty g^{-j}B$. Trivially, $B\subseteq\hat B$,
  and, by $\sigma$-additivity of the measure $\mu$, there is $s\in\nz$
  such that $\mu(\bigcup_{j=0}^sg^{-j}B)>\mu(\hat B)-\mu(B)$. As $\mu$
  is invariant under $g$ this implies, for all $r\in\nz$,
\begin{displaymath}
   \mu(B)+\mu\left(\bigcup_{j=0}^sg^{-r-j}B\right)
   =
   \mu(B)+\mu\left(\bigcup_{j=0}^sg^{-j}B\right)
   >
   \mu(\hat B)
\end{displaymath}
But $B\cup\bigcup_{j=0}^sg^{-r-j}B\subseteq\hat B$ so that
$B\cap\bigcup_{j=0}^sg^{-r-j}B\neq\emptyset$. We conclude that for each
$r\in\nz$ there is $j_r\in\{0,\dots,s\}$ such that $B\cap
g^{-r-j_r}B\neq\emptyset$ and hence $B\cap g^{r+j_r}B\neq\emptyset$. So
$R_\delta(x)$ has gaps bounded by $s$.
\end{proof}

Combining the last two lemmas with the elementary dichotomy
\eqref{eq:dichotomy} we arrive at the following conclusion which is
essentially Proposition 1 in \cite{GW1}.
\begin{prop}[\cite{GW1}]\label{prop:GlWe}
  Suppose that $X$ has no isolated points and that $g$ is transitive and has a
  finite invariant measure with full topological support. Then either
\begin{itemize}
\item[-] $g$ has sensitive dependence; or
\item[-] the family $(g^n)_{n\in\nz}$ is equicontinuous and $g:X\to X$
is a minimal homeomorphism.
In this case also the family $(g^n)_{n\in\zz}$ is equicontinuous.
\end{itemize}
\end{prop}
Observe that in the second case $g$ is uniquely ergodic.\footnote{
  This is well known and follows from the fact that ergodic limits of
  continuous functions are
  continuous (by equicontinuity) and hence constant (by transitivity).}\\ \\
\begin{proof}
  Let $x\in X$ be a transitive (and hence nonwandering) point and
  suppose that $g$ does not have sensitive dependence. In view of the
  dichotomy \eqref{eq:dichotomy}, $x$ is a Lyapunov point. Therefore
  the equicontinuity follows from Lemmas~\ref{lem:g2}
  and~\ref{lem:g3}. $g$ is minimal, \ie $\Tr=X$, because $X=\LS$ by
  equicontinuity and $\LS=\Tr$ by \eqref{eq:dichotomy}. Further, $g$
  is uniformly rigid by Lemma~\ref{lem:uniformlyrigid}, such that it
  must be a homeomorphism. Finally, the equicontinuity of
  $(g^n)_{n\in\BZ}$ follows again from the uniform rigidity of $g$:
  The iterates of any two points cannot come arbitrarily close to each
  other, as they become seperated again when $g^n$ is sufficiently
  close to the identity. But this implies the equicontinuity of the
  backwards iterates.
\end{proof}

\section{Sensitivity on SNAs}
\label{sec:sensitivity-SNA}

Now let $X=\tz\times M$, $M$ a metric space, and let $f:X\to X$ be a
continuous quasiperiodically forced map.  We assume throughout this
section that
\begin{eqnarray}
&&\triangleright\quad\label{property:1}
\text{$\A$ is a compact $f$-invariant subset of $X$ which has no
   isolated points,\hspace*{3.3cm}}\\
&&\triangleright\quad\label{property:2}
\text{that $f|_\A$ is transitive, and}\\
&&\triangleright\quad\label{property:3}
\text {that $\A$ is the topological support of a finite
$f$-invariant measure.}
\end{eqnarray}

\begin{rmk} \label{rmk:sna} 
An \emph{invariant graph} is usually defined as a measurable function
$\varphi:\tz \rightarrow M$ that satisfies
\[
f_\theta(\varphi(\theta)) \ = \ \varphi(\theta+\omega) \ \ \ \ \  \forall
\theta \in \tz  \ ,
\] 
with Lyapunov exponent $\lambda(\varphi)=\int_{\tz}
\log |f'_\theta(\varphi(\theta))| \ d\theta$ in case the fiber maps are
one-dimensional and differentiable. 

However, there is a subtle issue concerning this defintion which we do
not want to treat systematically here (this is done e.g.\ in \cite{J1}
at the end of section two), but nevertheless feel obliged to mention:
We do not want to distinguish between two invariant graphs which
coincide Lebesgue-a.e., and in particular we do not want to call an
invariant graph non-continuous if it is Lebesgue-a.s.\ equal to a
continuous one. Thus we will implicitely consider an invariant graph
to be an equivalence class of Lebesgue-a.s.\ equal graphs, and by the
topological closure of such an equivalence class we mean the smallest
compact set that can be obtained as the topological closure of a
representitive. This set coincides with the topological support of the
measure $\mu_\varphi$, which is obtained by projecting the Lebesgue
measure onto the invariant graph (here it does not matter which
representative is chosen). Note that $\mu_\varphi$ is ergodic w.r.t.\
$f$. 

If we now call the topological closure (in the above sense) of an
invariant graph with negative Lyapunov exponent a SNA, then it becomes
clear that this satisfies the assumptions
(\ref{property:1})--~(\ref{property:3}).
\end{rmk}
For any set $\A \subseteq X$, we denote by $\A_\theta$ its intersection with
the $\theta$-fiber, more precisely %
$\A_\theta := \{ \xi \in M \mid (\theta,\xi) \in \A\}$.  The following
concept turned out to be very important in the study of
quasiperiodically forced maps (see: \cite{stark:2003}):
\begin{definition}
  $\A \subseteq X$ is called \emph{pinched}, if for some $\theta \in
  \BT$ the set $\A_\theta$ consists only of a single point. In this case we
  call $\A$ pinched at $\theta$.
\end{definition}
Obviously, if $\A$ is invariant and pinched, then it is pinched on a
whole dense set, namely on the forward orbit of a pinched
fiber. If in addition $\A$ is compact then the set of $\theta$ at
which $\A$ is pinched is even residual. This follows quite easily from
a Baire argument, as in this case all sets $B_n := \{ \theta \in \BT
\mid \diam(\A_\theta) < \frac{1}{n} \}$ are open and dense, and their
intersection gives exactly the set of $\theta$ where $\A$ is pinched.

The next two results follow from the more general results of section two.
\begin{coro}\label{coro:1}
  Suppose $\A$ satisfies assumptions
  \eqref{property:1}\,--\,\eqref{property:3}. If $\A$ is pinched, then either
\begin{itemize}
\item[-] $f|_\A$ has sensitive dependence; or
\item[-] $\A$ is the graph of a continuous function.
\end{itemize}
\end{coro}
\begin{proof}
If $f|_\A$ does not have sensitive dependence, then $(f^n|_\A)_{n\in\nz}$ is
equicontinuous.  Let $\epsilon>0$. There is $\delta=\delta(\epsilon)>0$ such
that for all fibers $\A_\theta$ of diameter less than $\delta$ all their
images $f^n\A_\theta=\A_{\theta+n\omega}$ have diameter less than $\epsilon$.
As the set of all $\theta$ with $\diam(\A_\theta)<\delta$ is open and
nonempty, the minimality of the rotation $\theta\mapsto \theta+\omega$
implies that all fibers have diameter less than $\epsilon$. As $\epsilon>0$
was arbitrary, it follows that $\A$ is the graph of a function $\psi:\tz\to
M$. As $\A$ is compact, the function $\psi$ is continuous.
\end{proof}

\noindent
Sometimes the case where $\A$ is not necessarily pinched can also be
dealt with easily. Here is an example.
\begin{coro}\label{coro:2}
Suppose that $M$ is a compact interval and that all fiber maps are
monotone increasing and let $\A$ be as before.  Then either
\begin{itemize}
\item[-] $f|_\A$ has sensitive dependence; or
\item[-] $\A$ is the graph of a continuous function.
\end{itemize}
\end{coro}
\begin{proof}
Suppose that $f|_\A$ does not have sensitive dependence so that
$(f^n|_\A)_{n\in\nz}$ is equicontinuous and minimal by
Corollary~\ref{coro:1}.  Let $\A^+:=\{(\theta,\sup
\A_\theta):\theta\in\tz\}$, $\A^-:=\{(\theta,\inf
\A_\theta):\theta\in\tz\}$. As $f$ has monotone fiber maps, both,
$\A^+$ and $\A^-$ are $f$-invariant subsets of $X$. Hence, by
minimality of $f$, $\overline{\A^+}=X=\overline{\A^-}$. But this
implies that $\A$ is pinched\footnote{See \cite[Lemma
   4.3(i)]{FJJK-05} for the elementary proof.} so that
Corollary~\ref{coro:1} applies again.
\end{proof}

\noindent
Also in more delicate situations the dichotomy of Proposition~\ref{prop:GlWe}
can be useful.
\begin{prop}\label{prop:circle}
Suppose that $M=\tz$ and that all fiber maps are orientation preserving
circle homeomorphisms. Let $\A$ be a proper subset of $X$ which satisfies
assumptions \eqref{property:1}\,--\,\eqref{property:3}. Then either
\begin{itemize}
\item[-] $f|_\A$ has sensitive dependence; or
\item[-] $\A$ is the disjoint union of a finite number of disjoint
   curves which are cyclically permuted by the action of $f$.
\end{itemize}
\end{prop}
Information on situations where $\A=X$ will be provided in
section~\ref{sec:whole-space}.\\[2mm]
\begin{proof}
  Suppose that $f|_\A$ does not have sensitive dependence. Then the
  family $(f^n|_\A)_{n\in\BZ}$ is equicontinuous and $f|_\A$ is a
  minimal homeomorphism by Proposition~\ref{prop:GlWe}.

We introduce some more notation: Let $f^n_\theta(\xi) := \pi_2 \circ
f^n(\theta,\xi)$.  Further, for $\theta\in\BT$ let $\J_\theta$ be the family
of all connected components of $\BT\setminus \A_\theta$. So each
$J\in\J_\theta$ is a maximal interval in the complement of
$\A_\theta$.  Note that $f_\theta^nJ\in\J_{\theta+n\omega}$ if and
only if $J\in\J_\theta$.

For $\theta\in\BT$ and $J\in\J_\theta$ let
$s(J):=\sup_{n\in\zz}|f_\theta^nJ|$.
As the endpoints of such intervals $J$ belong to $\A$ and as the
family $(f^n|_\A)_{n\in\BZ}$ is equicontinuous, there is an
increasing function $\delta:(0,1]\to(0,1]$ such that
$|J|\geq\delta(t)$ whenever $s(J)\geq t > 0$.

For $t\geq0$ let
\begin{displaymath}
   N_t(\theta):=
   \begin{cases}
     \card\{J\in\J_\theta:s(J)\geq t\}&\text{if }t>0\\
     +\infty&\text{if }t=0\;.
   \end{cases}
\end{displaymath}
Note that for any $t >0$ there holds $0\leq N_t\leq\delta(t)^{-1}$
and that, for each fixed $\theta$, $t\mapsto N_t(\theta)$ is a
decreasing function continuous from the left. As $s(f_\theta
J)=s(J)$ for all $J\in\J_\theta$, we see that
$N_t(\theta+\omega)=N_t(\theta)$.

Next, for $p=1,2,3,\dots$, let
\begin{displaymath}
   \gamma_p(\theta):=\sup\{t\in\BR:N_t(\theta)\geq p\}\;.
\end{displaymath}
Obviously $0\leq\gamma_p\leq1$ and
$\gamma_p(\theta+\omega)=\gamma_p(\theta)$. Further, as $t\mapsto N_t(\theta)$
is continuous from the left we have $N_{\gamma_p(\theta)}(\theta)
\geq p$.

The function $\gamma_1:\BT\to\BR$ plays a special role: Observe
first that
\begin{displaymath}
   \gamma_1(\theta)=\sup_{n\in\BZ}\ell(\theta+n\omega)\text{\; where\; }
   \ell(\theta):=\max\{|J|:J\in\J_\theta\}\;.
\end{displaymath}
As $\A\subseteq X$ is closed, the function $\ell:\BT\to\rz$ is lower
semicontinuous and so are the functions
$\theta\mapsto\ell(\theta+n\omega)$.  Hence, as a supremum of lower
semicontinuous functions, also $\gamma_1$ is lower semicontinuous,
and as $\gamma_1$ is invariant under rotation by the irrational
$\omega$, it must be constant.  As $\A\neq X$ by assumption, we have
$\gamma_1>0$.  We turn to the other functions $\gamma_p:\BT\to[0,1]$.
Let
$\delta_0:=\delta(\frac {\gamma_1}2)\leq\frac{\gamma_1}2$.\\[2mm]
\textbf{Claim:} The sets $\{\gamma_p\leq c\}$ are
closed for all $c\in(\gamma_1-\delta_0,\gamma_1]$.\\[2mm]
Indeed, consider a sequence of $\theta_k\in\{\gamma_p\leq c\}$ which
converge to some $\theta\in\BT$. Let $t\in(c,\gamma_1]$ and denote
$q:=N_t(\theta)$. Then $q\geq1$ and there are pairwise disjoint
$J_1,\dots,J_{q}\in\J_\theta$ with $s(J_i)\geq t$ for all $i$.
Consider compact subintervals $\hat J_i\subset J_i$. As $X\setminus
\A$ is open there are, for all sufficiently large $k\in\nz$,
intervals $J_1^k,\dots,J_{q}^k\in\J_{\theta_k}$ such that $\hat
J_i\subset J_i^k$ for $i=1,\dots,q$.  Since this holds for all
choices of the compact subintervals $\hat J_i$, we conclude
\begin{displaymath}
   \liminf_{k\to\infty}|f_\theta^nJ_i^k|\geq|f_\theta^nJ_i|\quad\text{ for
   all $i=1,\dots,q$ and all $n\in\nz$}
\end{displaymath}
and therefore
\begin{displaymath}
   \liminf_{k\to\infty}s(J_i^k)\geq s(J_i)\geq t>c\quad\text{ for
   all $i=1,\dots,q$.}
\end{displaymath}
This allows the conclusion
$N_t(\theta)=q\leq\liminf_{k\to\infty}N_u(\theta_k)\leq p-1$ for all
$t>u>c$ so that $\gamma_p(\theta)\leq c$
-- and thus proves the claim -- once we have shown that the intervals
$J_1^k,\dots,J_{q}^k$ are pairwise distinct for large $k$. Suppose
for a contradiction that this is not the case. Then, w.l.o.g.,
$J_1^k=J_2^k$ for infinitely many $k$. Hence $\hat J_1\cup\hat
J_2\subseteq J_1^k$ for infinitely many $k$ and so
\begin{displaymath}
   \gamma_1
   \geq
   \limsup_{k\to\infty}|f_\theta^nJ_1^k|
   \geq
   |f_\theta^n\hat J_1|+|f_\theta^n\hat J_2|\quad\text{ for all $n\in\zz$.}
\end{displaymath}
Again this holds for all choices of the compact subintervals $\hat
J_1,\hat J_2$, so that
\begin{displaymath}
   \gamma_1
   \geq
   |f_\theta^n J_1|+|f_\theta^n J_2|\quad\text{ for all $n\in\zz$.}
\end{displaymath}
But $s(f_\theta^nJ_2)=s(J_2)\geq
t>\gamma_1-\delta_0\geq\frac{\gamma_1}2$, whence
$|f_\theta^nJ_2|\geq\delta(\frac{\gamma_1}2)=\delta_0$. This yields
the contradiction
\begin{displaymath}
   \gamma_1
   \geq
   s(J_1)+\delta_0
   \geq
   t+\delta_0
   >
   c+\delta_0>\gamma_1
\end{displaymath}
and finishes the proof of the claim.

We can summarize that for $c\in(\gamma_1-\delta_0,\gamma_1]$ the
closed sets $\{\gamma_p\leq c\}$ are invariant under rotation by the
irrational $\omega$, so for all these $c$ and for all
$p\in\{1,2,3,\dots\}$, either $\{\gamma_p\leq c\}=\emptyset$ or
$\{\gamma_p\leq c\}=\BT$.  Let
$q:=\max\{p\in\BN:\exists\theta\in\BT\text{ s.t.
}\gamma_p(\theta)=\gamma_1\}$. Then $\gamma_q=\gamma_1$ is constant
and there is $\eta>0$ such that $\gamma_{q+1}\leq\gamma_q-\eta$.
This means that for each $\theta$ there are $q$ intervals
$J\in\J_\theta$ with $s(J)=\gamma_1$ and the $s$-value of all other
intervals is at most $\gamma_1-\eta$.

Now we are ready to finish the proof of the proposition. Let $E_\theta^\pm$
be the sets of the $q$ ``upper'' respectively ``lower'' endpoints of those
intervals $J\in\J_\theta$ with $s(J)=\gamma_1$. We show that the set valued
maps $\theta\mapsto E_\theta^+$ and $\theta\mapsto E_\theta^-$ are
continuous: As above consider a sequence of $\theta_k\in\BT$ which converge
to some $\theta\in\BT$.  There are pairwise disjoint intervals
$J_1,\dots,J_{q}\in\J_\theta$ with $s(J_i)=\gamma_1$ for all $i$.  Consider
compact subintervals $\hat J_i\subset J_i$. As $X\setminus \A$ is open there
are, for all sufficiently large $k\in\nz$, intervals
$J_1^k,\dots,J_{q}^k\in\J_{\theta_k}$ such that $\hat J_i\subset J_i^k$ for
$i=1,\dots,q$. Now let $\epsilon\in(0,\eta)$ and choose $\delta\in(0,\eta)$
as in the definition of equicontinuity \eqref{eq:equicontinuous-def}.  As
$s(J_i)=\gamma_1$ for all $i$, there are $n_1,\dots,n_q\in\BZ$ such that
$|f^{n_i}J_i|\geq \gamma_1-\frac\delta2$, while
$|f^{n_i}J_i|,|f^{n_i}J_i^k|\leq \gamma_1$ for all $i$ and $k$ by definition
of $\gamma_1$. We can choose the intervals $\hat J_i\subset J_i$ such that
$|f^{n_i}\hat J_i|> \gamma_1-\delta$. So also
$|f^{n_i}J_i^k|>\gamma_1-\delta$ for large $k$.  This has two implications
which together yield the continuity of the set valued maps:
$s(J_i^k)>\gamma_1-\eta$ and hence $s(J_1^k)=\dots=s(J_q^k)=\gamma_1$ for
large $k$, and second, the corresponding endpoints of the intervals
$f^{n_i}J_i$ and $f^{n_i}J_i^k$ have distance less than $\delta$, so that
the corresponding endpoints of the intervals $J_i$ and $J_i^k$ have distance
at most $\epsilon$.

The graphs of the maps $\theta\mapsto E_\theta^+$ and $\theta\mapsto
E_\theta^-$ are thus closed invariant subsets of $\A$ so that, by
minimality of $f|_\A$, both graphs are identical and coincide with
$\A$.  It follows that $\card(\A_\theta)=\gamma_1$ and
$\inf\{d(\xi,\zeta): \xi,\zeta\in \A_\theta,
\xi\neq\zeta\}\geq\delta(\gamma_1)>0$. From this the second
alternative of the proposition follows.
\end{proof}

\begin{rmk}
  The results of this section remain valid if the rotation
  $\theta\mapsto\theta+\omega$ which forces the system is replaced by
  any transitive continuous map $R$ on a compact metric space $Z$
  which has no isolated points and which admits a finite invariant
  measure with full topological support. Indeed, if the forced system
  does not have \sdic, then also $R:Z\to Z$ does not have \sdic\ so
  that $R$ is a minimal homeomorphism of $Z$ and $(R^n)_{n\in\zz}$ is
  equicontinuous, see Proposition~\ref{prop:GlWe}. But only these two
  properties of the rotation, minimality and equicontinuity, were used
  in the proofs of this section, so the proofs carry over without
  changes to the more general $R$.
\end{rmk}
\section{Sensitivity on the whole phase space}\label{sec:whole-space}

In this section we turn to the question of sensitive dependence on the
whole phase space. In order to do so, we restrict to two classes of
quasiperiodically forced systems, namely quasiperiodically forced
circle homeomorphisms and quasiperiodically forced monotone interval
map.

As in the last section, we will use the notation $f^n_\theta(\xi)
:= \pi_2 \circ f^n(\theta,\xi)$.  We say $f$ is a
\textit{quasiperiodically forced circle homeomorphism} if $M=\BT1$
and each fiber map $f_\theta$ is a homeomorphism of the circle. By
$F$ we will denote a continuous lift of $f$ to $\BT1\times\BR$. First
of all, the case where $f$ is not homotopic to the identity can be
treated quite easily:
\begin{prop} \label{lem:nonhomotopic}
Suppose $f$ is a quasiperiodically forced circle homeomorphism which
is not homotopic to the identity. Then $f$ has \sdic\ on $\BT2$.
\end{prop}
\proof \\
As $f$ is not homotopic to the identity, it is transitive (see
\cite{kim/kim/hunt/ott:2003}), such that we can apply
Lemma~\ref{lem:uniformlyrigid} to see that either $f$ has \sdic\ on
$\BT2$ or $f$ is uniformly rigid. Suppose $f$ is uniformly rigid.
Then $f^n$ is arbitrarily close to the identity for infinitly many $n
\in \BN$, and in particular the image of a constant line $\Gamma =
\BT1\times\{\xi\}$ is mapped arbitrarily close to itself by $f^n$.
However, this would imply that $\Gamma$ and $f^n(\Gamma)$ are in the
same homotopy class, contradicting the fact that $f$ is not homotopic
to the identity. Therefore $f$ must have \sdic.

\qed

\noindent
Now suppose $f$ is homotopic to the identity. In this case, Herman
showed in \cite{herman:1983} that similar to the unforced case
the limit
\begin{equation}
   \label{eq:rotnum}
       \rho_{F} := \lim_{n\rightarrow \infty}
       \frac{1}{n}(F_{\theta}^n(\xi) - \xi)
\end{equation}
exists for any continuous lift $F$ of $f$ and is independent of
$\theta$ and $\xi$. Further, $\rho_f := \rho_F \bmod 1$ does not
depend on the choice of the lift $F$. However, unlike unforced
circle homeomorphisms the so-called \textit{deviations from the
  constant rotation}
\begin{equation} \label{eq:deviations}
     |F_\theta^n(\xi) - \xi - n\rho_F|  \ ,
\end{equation}
need not be bounded uniformly in $\theta, \xi$ and $n$, and in fact
an important distinction can be made with respect to this: $f$ is
called $\rho$-bounded if the quantities in (\ref{eq:deviations}) are
uniformly bounded and $\rho$-unbounded otherwise. If the systems is
$\rho$-bounded, the dynamics can be understood quite easily: In this
case an analogue to Poincar\'e's famous classification of the dynamics
of circle homeomorphisms holds, such that the system is either
semi-conjugate to an irrational translation of the torus and $\rho_f$ is not
rationally related to the rotation number $\omega$ on the base, or
there exists an invariant strip, which is the suitable
analogue for a fixed or periodic point in this setting (see
\cite{jaeger/stark:2005}), and the rotation numbers $\rho_f$ and
$\omega$ are rationally related.

The more interesting case, which does not occur in the one-dimensional
situation and which we will consider in the following, is the
$\rho$-unbounded one. Here neither of the two above alternatives can
occur, the system is always topologically transitive (see
\cite{jaeger/stark:2005}), and as we will see below it also has \sdic\
on the whole phase space. However, before we can show this we need the
following statement, which is contained in
\cite{stark/feudel/glendinning/pikovsky:2002}:
\begin{lem}
  \label{lem:sfgp} Suppose $F$ is the lift of a quasiperiodically forced circle
  homeomorphism homotopic to the identity which is
  $\rho$-unbounded. Then there exists a residual set of $\theta$
  such that the deviations
   \begin{equation} \label{eq:deviationsII}
     F_\theta^n(\xi) - \xi - n\rho_F  \
\end{equation}
are unbounded both above and below (independent of $\xi$), but at the
same time there exist two disjoint dense sets of $\theta$ such
that the deviations (\ref{eq:deviationsII}) are bounded uniformly
from above, respectively below. (However, there exists no orbit
on which the deviations are bounded both above and below.)
\end{lem}
Now we can prove the following:
\begin{prop} \label{lem:rhounbounded}
Suppose $f$ is a quasiperiodically forced circle homeomorphism,
homotopic to the identity, which is $\rho$-unbounded. Then $f$ has
\sdic\ on $\BT2$.
\end{prop}
\proof \\
Let $F : \BT1 \times \BR\to\BT1 \times \BR$ be a lift of $f$ and choose
$\epsilon \in (0,\frac{1}{4})$ such that $d(\xtil,\ytil) < \epsilon$
implies $d(F(\xtil),F(\ytil)) < \frac{1}{4}$ for all $\xtil,\ytil \in
\BT1 \times \BR$.  Note that $d(\xtil,\ytil) < \frac{1}{4}$ implies
$d(\pi(\xtil),\pi(\ytil)) = d(\xtil,\ytil)$. We will now show that $f$ is
$\epsilon$-sensitive on $\BT2$, that is $\SD_\epsilon = \BT2$.

To that end, choose any $x \in\BT2$ and $\delta>0$. Let $\xtil \in
\BT1 \times \BR$ be a lift of $x$, i.e.\ $\pi(\xtil) = x$. As both
the fibers which are $\rho$-bounded above and those which are
$\rho$-unbounded above are dense, we can find both a point $\ytil$
which is $\rho$-bounded above and a point $\ztil$ which is
$\rho$-unbounded above in $B_\delta(\xtil)$. As $\sup_{n\in\BN}
d(F^n\ytil,F^n\ztil) = \infty$ and due to the choice of $\epsilon$,
this means that for some $m \in \BN$ we must have
$d(F^m\ytil,F^m\ztil) \in [\epsilon,\frac{1}{4})$.  Let $y=\pi(\ytil)$
and $z=\pi(\ztil)$. Then $d(F^my,F^mz) \in [\epsilon,\frac{1}{4})$ as
well, and as $y,z \in B_\delta(x)$ this completes the proof.

\qed

\noindent It is this proposition, which applies to Herman's examples
mentioned in the introduction. In \cite{herman:1983} Herman studies
$\textrm{SL}(2,\BR)$-cocycles over irrational rotations, that is
mappings $(\om,A) :\BT \times \BR2 \rightarrow \BT \times \BR2$,
$(\theta,v) \mapsto (\theta+\om, A(\theta)v)$ where $\om \in \BT$ is
irrational and $A : \BT \rightarrow \textrm{SL}(2,\BR)$ is a
continuous matrix-valued function. By their action on the real
projective space and subsequent identification of $\mathbb{P}(\BR2)$
with $\BT$, such a cocycle $(\omega,A)$ induces a quasiperiodically
forced circle homeomorphism $f_A$. If the Lyapunov exponent
$\lambda(\omega,A) := \liminf_{n\rightarrow\infty} \frac{1}{n}
\int_{\BT} \log \| A(\theta+(n-1)\omega) \circ \ldots \circ A(\theta)
\| \ d\theta$ of such a cocycle is positive, the corresponding cocycle
$f_A$ will have exactly two invariant graphs, one with positive and
one with negative Lyapunov exponent (corresponding to the stable and
unstable subspaces in Oseledets Multiplicative Ergodic Theorem).
Herman showed that for any pair of rotation numbers $\omega$ and
$\rho$ there exist cocycles $(\omega,A)$ such that $f_A$ has fiberwise
rotation number $\rho$ and the Lyapunov exponent of the cocycle is
positive (Proposition 4.6 in \cite{herman:1983}). If the rotations
numbers are chosen rationally independent, the corresponding map $f_A$
must be $\rho$-unbounded: As there exist invariant graphs it cannot
be semi-conjugate to an irrational torus translation, and as the
rotation numbers are not rationally dependend there cannot be any
invariant strips. Thus, both alternatives in the $\rho$-bounded case
are ruled out, so these examples are
$\rho$-unbounded, topologically transitive, and by the preceding
proposition they have \sdic\ on $X$. The only candidates for proper
minimal subsets of the whole space are the essential closures of the two
invariant graphs%
\footnote{The topological support of the measures which are obtained by
  projecting the Lebesgue measure on the base onto the invariant graphs.}
and restricted to these $f_A$ has \sdic\ as well by
Proposition~\ref{prop:circle}. However, it should be mentioned that it
is still an open question whether transitive but non-minimal dynamics
do really exist in this setting. To the knowledge of the authors, the
only examples where the topological dynamics of such cocycles has been
clarified so far are certain quasiperiodic Schr\"odinger cocycles, for
which Bjerkl\"ov proved under some additional assumptions that the
dynamics are minimal \cite{bjerkloev:2005a}.

We now consider the case in which the fiber maps $f_\theta$ are maps
of the interval. In this case, continuous and non-continuous invariant
graphs may coexist, and consequently there might be different regions
in the phase space with and without \sdic. Hence, instead of looking
at the whole phase space we concentrate on the `domain of attraction'
of a non-continuous invariant graph.
\begin{lem} \label{lem:domainofattraction}
Suppose $f$ is a quasiperiodically forced monotone interval map and
$\varphi$ is an upper semi-continuous, non-continuous invariant
graph. Then there exists $\epsilon >0$ such that the following
inclusion holds:

\[
A_\varphi \ := \ \{(\theta,\xi) \mid \xi > \varphi(\theta)
\textrm{ and } \inf_{n\in\BN} |f_\theta^n(\xi) -
\varphi(\theta+n\omega)| = 0 \} \ \subseteq \ \SD_\epsilon \ .
\]
\end{lem}
\proof \\
Let $\Phi := \{(\theta,\varphi(\theta) \mid \theta \in \BT\}$ and
denote its topological closure by $\overline{\Phi}$.  Let
$\varphi^-(\theta) := \inf\{\xi\mid(\theta,\xi) \in
\overline{\Phi}\}$. Then $\varphi^-$ is a lower semi-continuous
invariant graph and the set
$[\varphi^-,\varphi] := \{ (\theta,\xi) \mid \varphi^-(\theta)
\leq \xi \leq \varphi(\theta)]$ is pinched, i.e.\
$\varphi^-(\theta) =
\varphi(\theta)$ for a residual set of $\theta$.%
\footnote{This is quite easy to see using a Baire argument, see
  \cite[Lemma 4.3(i)]{FJJK-05} or \cite{stark:2003} for details.}
Choose some $\theta_0 \in \BT1$ with $\varphi^-(\theta_0) \neq
\varphi(\theta_0)$ and let %
$\epsilon := \frac{1}{4} (\varphi(\theta_0) -
\varphi^-(\theta_0))$. Further, let
$x = (\theta,\xi) \in A_\varphi$ and $\delta > 0$ be given. As
$\varphi$ is upper semi-continuous, we can assume w.l.o.g\ (by
decreasing $\delta$ if necessary) that $\varphi(\tilde{\theta}) <
\xi \ \forall \tilde{\theta} \in B_\delta(\theta)$.

Due to the definition of $\varphi^-$, we can find $\theta_1$ in
$B_{\delta/2}(\theta_0)$ with $\varphi(\theta_1) -
\varphi^-(\theta_0) \leq \epsilon$, such that $\varphi(\theta_1)
\ \leq \ \varphi(\theta_0) - 3 \epsilon$. Further, as $\varphi$ is
upper semi-continuous there exists some $\eta \in (0,\delta/2)$
such that
\begin{equation}
\label{eq:doa1}
\varphi(\tilde{\theta}) \ \leq \ \varphi(\theta_1) + \epsilon \
\leq \ \varphi(\theta_0) - 2
\epsilon \ \ \ \forall \tilde{\theta} \in B_{\eta}(\theta_1) \ .
\end{equation}
Now we choose some $n\in\BN$ which satisfies $\theta + n\omega \in
B_{\eta}(\theta_1)$ and
\begin{equation}
\label{eq:doa2}
|f_\theta^n(\xi)-\varphi(\theta+n\omega)| \ \leq \ \epsilon \ .
\end{equation}
Such an integer exists because the set $\{ k\in\BN \mid
\theta+k\omega\in B_{\eta}(\theta_1) \}$ has bounded gaps (in the
sense of (\ref{eq:bounded-gaps})) and the orbit of $x\in A_\varphi$
will stay $\epsilon$-close to $\varphi$ for arbitrarily long time
intervals due to the definition of $A_\varphi$ and the continuity of
$f$. Consequently we obtain
\begin{equation}
\label{eq:doa3}
f_\theta^n(\xi) \ \leq \ \varphi(\theta+n\omega)+\epsilon \ \leq \
\varphi(\theta_0) - \epsilon
\end{equation}
by (\ref{eq:doa2}) and (\ref{eq:doa1}). At the same time $
f^n_{\theta_0-n\omega}(\xi) \geq \varphi(\theta_0)$, such that
$|f^n_\theta(\xi)-f^n_{\theta_0-n\omega}(\xi)| \geq \epsilon$.  As
$y=(\theta_0-n\omega,\xi) \in B_\delta(x)$ (note that $\theta+n\omega
\in B_{\eta}(\theta_1) \subseteq B_\delta(\theta_0)$, such that
$\theta_0-n\omega \in B_\delta(\theta)$), this completes the proof.

\qed

\noindent
Obviously, an analogous statement holds for the region below a lower
semi-continuous invariant graph. As an application we obtain the
following proposition, which in particular contains the second
statement about pinched systems mentioned in the introduction.
\begin{cor}
\label{cor:pinchedga}
Suppose $f$ is a quasiperiodically forced monotone interval map, such
that the global attractor ${\cal K}:= \bigcap_{n\in\BN}
f^n(\BT1\times[a,b])$ is pinched and the upper and lower bounding
graphs $\varphi^+(\theta) := \sup {\cal K}_\theta$ and $\varphi^-(\theta) :=
\inf {\cal K}_\theta$ are non-continuous. Then $f$ has \sdic\ on
$\BT1\times[a,b]$. The same is true if one of the bounding graphs is
continuous, but coincides with one of the boundaries of the annulus.
\end{cor}
\proof \\
We treat the case of two non-continuous bounding graphs, the second
case is similar. Any point $(\theta,\xi)$ above $\varphi^+$ is
necessarily contained in $A_{\varphi^+}$. Thus we can apply the above
Lemma~\ref{lem:domainofattraction} to see that for some suitable
$\epsilon > 0$ we have $\{ (\theta,\xi) \mid \xi > \varphi^+(\theta)\}
\ \subseteq \ \SD_\epsilon$.  Similarly, we can assume $\{\theta,\xi)
\mid \xi < \varphi^-(\theta)\} \ \subseteq \ \SD_\epsilon$, such that
together we have ${\cal K}^c \subseteq \SD_\epsilon$. But as ${\cal K}$ is
pinched and therefore has empty interior and $\SD_\epsilon$ is closed,
this implies $\BT1\times[a,b] \subseteq \SD_\epsilon$.

\qed


\section{A final example}
\label{sec:final}
Given the great attention pinched skew products have received after
they had been introduced by Grebogi et al.\ in \cite{GOPY}, it is
rather surprising that the second type of model system which was
proposed in the very same paper has been completely neglected so far.
As Example 2 in \cite{GOPY} the authors consider the map $\Lambda : \BT1
\times \BR2 \rightarrow \BT1 \times \BR2$ depending on parameters
$\beta$ and $\gamma$ and given by
\begin{equation}
\label{eq:Lambda}
\Lambda(\theta,\xi) \ := \ \left( \theta + \omega,
   \frac{\beta}{1+\twonorm{\xi}^2} \cdot  \Matrix{cc}{1 & 0 \\ 0 &
     \gamma} \cdot  R_\theta\cdot  \xi \right) \ .
\end{equation}
Here $R_\theta$ denotes the rotation matrix
\[
R_\theta \ := \
\Matrix{cc}{\cos(2\pi\theta) & \sin(2\pi\theta)\\
-\sin(2\pi\theta) & \cos(2\pi\theta)} \
\]
and $\xi={u\choose v}$ is a vector as in Figure~\ref{fig:1}.
In order to obtain a compact phase space, we choose a sufficiently
large constant $C>0$ such that $X:= \BT1 \times B_C(0)$ is mapped
strictly inside itself (i.e.\ $\Lambda(X) \subseteq \textrm{int}(X)$)
and consider $\Lambda$ restricted to $X$.

Similar to pinched skew products, the 0-line $\xi=0$ is invariant.
Further, as the action of $\Lambda$ on any continuous curve that does
not intersect the 0-line increases the number of lefthand turns around
the 0-line there can be no other continuous invariant curve (in other
words, the projective action of $\Lambda$ is not homotopic to the
identity). The numerical results in \cite{GOPY} indicate that for the
considered parameter values ($\beta=2,\ \gamma=0.5$ and $\omega$ the
golden mean) the system exhibits an SNA. This SNA seems to be a
quasiperiodic two-point attractor (i.e.\ a two-valued measurable
invariant graph) which attracts Lebesgue-a.e.\ initial condition. In
the following, we will give a rigorous proof of this observation and
show in addition that the attractor is embedded in a two-dimensional
torus ${\cal T}_0$, which is the boundary of the global attractor ${\cal
G} := \bigcap_{n\in\BN} \Lambda^n(X)$ in the three-dimensional phase
space. Further, $\Lambda$ has \sdic\ both restricted to the attractor
and on the whole phase space. These results remain valid as
long as $1<\beta\leq2,\, \gamma \in (0,1)$, and $\beta \gamma \geq 1$.

\ \\
\textit{A two-to-one factor.} In order to analyze the dynamics of
$\Lambda$, it turns out to be more convenient to use
polar coordinates in $\BR2 \setminus \{0\}$, and to
consider directions only projectively, rather than use the
standard Cartesian coordinates
in $\BR2$. Therefore, we will now
introduce a map $\tilde{\Lambda} : \BT2\times[0,C] \rightarrow
\BT2\times[0,C]$ which is a two-to-one factor of $\Lambda$.%
\footnote{To be absolutely precise, $\Lambdatil_{|\BT2 \times (0,C]}$
   will be a two-to-one factor of $\Lambda_{|\BT1 \times
     \BR2\setminus\{0\}}$, whereas the 0-line is `blown up' into the
     0-torus $S=\BT2 \times \{0\}$. However, as the 0-line is
     invariant and we are only  interested in the dynamics off the
     0-line, this is sufficient for our purposes.}
It will turn out that there exists an attracting invariant graph for
$\Lambdatil$, and the preimage of this graph under the factor map then
gives the two-point-attractor for $\Lambda$. However, we will have to
leave open here whether this attractor can further be decomposed into
two one-valued invariant graphs or not.

\ \\
Let $b(x) := \frac{x}{1+x2}$ and
\begin{equation} \label{eq:cocycle}
   A(\theta) \ := \
   \Matrix{cc}{\gamma^{-\frac{1}{2}}&0\\0&\gamma^{\frac{1}{2}}} \cdot
   R_\theta \ \in
   \ \textrm{SL}(2,\BR) \ .
\end{equation}
Then (\ref{eq:Lambda}) becomes
\begin{equation}
\label{eq:LambdaII}
\Lambda(\theta,\xi) \ = \
\left(\theta+\omega,\beta\gamma^{\frac{1}{2}} \cdot b(\twonorm{\xi})
   \cdot  A(\theta) \cdot \frac{\xi}{\|\xi\|_2} \right) \ .
\end{equation}
As mentioned, we will consider projective polar coordinates $\alpha =
\frac{1}{\pi}\cot^{-1}(\frac{u}{v}) \in \BT1$ and $r=\twonorm{\xi}$
for $\xi={u\choose v}\in\BR2 \setminus \{0\}$. The reason for doing
so is the fact that the action of $\Lambda$ on $\alpha$ does not
depend on $r$, such that the system becomes a skew product over a skew
product. Further, the dynamics of $\alpha$ are determined by the
projective action of the quasiperiodic $\textrm{SL}(2,\BR)$-cocycle
\begin{equation}
\label{eq:cocycleII}
(\omega,A) : \BT1 \times \BR2 \rightarrow \BT1 \times \BR2 \ \ \
, \ \ \ (\theta,\xi) \mapsto (\theta+\omega,A(\theta)\cdot \xi)
\ ,
\end{equation}
which induces  a quasiperiodically forced circle homeomorphisms
$f=f_A$.%
\footnote{If $A(\theta) = \left(\begin{array}{cc} a_\theta &
     b_\theta \\ c_\theta & d_\theta \end{array} \right)$,
then we can first define a map $\tilde{f}_A : \BT1 \times
\overline{\BR} \rightarrow \BT1 \times \overline{\BR}$ by
\[
    \tilde{f}_A(\theta,x) \ = \ \left( \theta + \omega,
      \frac{a_\theta x +  b_\theta}{c_\theta x +d_\theta}
    \right) \ ,
\]
and identification of $\overline{\BR}$ with $\BT1$  via $x
\mapsto \frac{1}{\pi}\cot^{-1}(x)$ yields $f=f_A$.  }
Such cocycles present one of the few classes of quasiperiodically
forced systems which are already well-understood, and in particular we
can apply results from \cite{herman:1983} and \cite{thieullen:1997} to
our problem.

If we now let $\Theta=(\theta,\alpha)$, we obtain a map $\Lambdatil :
\BT2 \times [0,C] =: Y \rightarrow Y$ given by
\begin{equation}
\label{eq:Lambdatil}
\tilde{\Lambda}(\Theta,r) \ := \ (f(\Theta),g_\Theta(r))
\end{equation}
where $g_\Theta$ is defined by the dynamics of $\Lambda$ on $r$. More
precisely, suppose $\xi={u\choose v} \in \BR2 \setminus \{0\}$ is a
vector with $\frac{1}{\pi}\cot^{-1}\left(\frac{u}{v}\right) =
\alpha$ and length $\twonorm{\xi}=1$ and let
\begin{equation}
\label{eq:adef}
a(\Theta) \ := \ \twonorm{\beta \gamma^{\frac{1}{2}} \cdot
     A(\theta)\cdot \xi} \ .
\end{equation}
Then it is easy to see from (\ref{eq:LambdaII}) that
\begin{equation}
\label{eq:guform}
g_\Theta(r) \ = \ a(\Theta) \cdot b(r) \
\end{equation}
and $a$ depends continuously on $\Theta$. Further, let
$X' := X \setminus (\BT1 \times \{0\})$ and
$Y' := Y \setminus (\BT2 \times \{0\})$. Then, as mentioned before,
$\Lambdatil_{|Y'}$ is a two-to-one factor of $\Lambda_{|X'}$ with
factor map
\begin{equation}
\label{eq:factormap}
h : (\theta,\xi) \mapsto \left(\theta,\frac{1}{\pi}
   \cot^{-1}\left(\frac{u}{v}\right),\twonorm{\xi}\right) \ .
\end{equation}

\ \\
\textit{Base dynamics of $\Lambdatil$.} First we analyze the
dynamics of the driving homeomorphism $f$. As we have already argued
above, $f$ is not homotopic to the identity and therefore
topologically transitive \cite{kim/kim/hunt/ott:2003} and has \sdic\
due to Lemma~\ref{lem:nonhomotopic}. Further, the Lyapunov exponent of
the cocycle $(\omega,A)$ is defined as
\begin{equation}
\label{eq:cocyclelyap}
\lambda(\omega,A) \ := \ \lim_{n\rightarrow \infty} \int_{\BT1}
   \log \parallel A_n(\theta) \parallel \ d\theta \ ,
\end{equation}
where $A_n(\theta) = A(\theta+(n-1)\omega) \circ \ldots \circ
A(\theta)$. Section 4.1 in \cite{herman:1983} provides a lower
bound for the Lyapunov exponent, namely
\begin{equation}
\label{eq:lyapbound}
\lambda(\omega,A) \ \geq \ \log\left(\frac{\sqrt{\gamma}}{2} +
   \frac{1}{2\sqrt{\gamma}}\right) \ .
\end{equation}
This means that for $\gamma \neq 1$ the Lyapunov exponent is always
positive, and Oseledet's Multiplicative Ergodic Theorem then implies
that there exists an invariant splitting of $\BR2$ into a stable and
an unstable subspace. This in turn is equivalent to the existence of
exactly two invariant graphs $\varphi^s$ and $\varphi^u$ for the
induced map $f$ with positive and negative Lyapunov exponent,
respectively.%
\footnote{In fact, we have $\lambda(\varphi^s) = 2\lambda(\omega,A)$
and $\lambda(\varphi^u) = -2\lambda(\omega,A)$ where $\varphi^s$ is
the invariant graph corresponding to the stable direction and
$\varphi^u$ the one corresponding to the unstable direction.}
Note that the invariant graph $\varphi^u$ corresponding to the
unstable subspace is the one with negative Lyapunov exponent and
attracts Lebesgue-a.e.\ initial condition.

\ \\
\textit{Ergodic invariant measures and vertical Lyapunov exponents.}
In order to obtain more information about our system, we have to
characterize the ergodic invariant measures for $\Lambdatil$. Further,
we have to determine their \textit{`radial'} Lyapunov exponents, which are
defined as
\begin{equation}
\label{eq:verticallyaps}
\lambda_{\textrm{rad}}(\nu) \ := \ \int_{Y} \log Dg_\Theta(r) \ d\nu(\Theta,r)
\ ,
\end{equation}
where $\nu$ is the invariant measure and $Dg_\Theta$ denotes the derivative
of $g_\Theta$ w.r.t.\ $r$.%
\footnote{We denote these Lyapunov exponents by
$\lambda_{\textrm{rad}}$ in order to distinguish them from the
Lyapunov exponents in two-dimensional skew products, thus avoiding
ambiguities.}
First of all, for the base dynamics given by $f$ there exist exactly
two ergodic invariant measures $\mu^s$ and $\mu^u$ which are
associated to the invariant graphs by
\begin{equation}
\label{eq:associatedmeasures}
\mu^i(A) \ := \ m(\{ \theta \in \BT1 \mid
(\theta,\varphi^i(\theta)) \in A \}) \ ,
\end{equation}
where $m$ denotes the Lebesgue measure on $\BT1$. These two measures
$\mu^i$ can be naturally identified with $\Lambdatil$-invariant
measures $\mutil^i$ by embedding them into the invariant 0-torus
$S:=\BT2 \times \{0\}$ in the canonical way. Their Lyapunov exponents are
then given by
\begin{equation}
\label{eq:muilyaps}
\lambda_{\textrm{rad}}(\mutil^i) \ = \ \int_{\BT1}
\log Dg_{(\theta,\varphi^i(\theta))}(0) \ d\theta \ = \
\int_{\BT1} \log a(\theta,\varphi^i(\theta)) \ d\theta \ .
\end{equation}
It is not hard to see that for the considered parameter
values $1<\beta\leq2$, $\beta\gamma \geq 1$, both exponents are positive:
{}From (\ref{eq:cocycle}) and (\ref{eq:adef}) we deduce that $1\leq a \leq 2$
and $a(\theta,\alpha)=1$ if and only if $\beta\gamma=1$ and $\alpha =
\frac{1}{2}-\theta$. As the invariant graphs are non-continuous and
therefore $\varphi^i(\theta) = \frac{1}{2} - \theta$ cannot hold
$m$-a.s., this implies that $\lambda_{\textrm{rad}}(\mutil^i) > 0$.

Any other ergodic invariant measure $\nu$ must project down to an
ergodic measure for the base dynamics, that is either $\mu^s$ or
$\mu^u$. Thus, in order to study $\nu$ we can restrict the base
dynamics to the respective invariant graph $\Phi^i$. But this means we
obtain a system which can be viewed as a two-dimensinal skew product
$h^{(i)} : \BT1 \times [0,C]$ again, with fiber maps
\begin{equation}
\label{eq:hfibers}
h^{(i)}_\theta(r) \ = \ g_{(\theta,\varphi^i(\theta))}(r) \ = \
a(\theta,\varphi^i(\theta)) \cdot b(r) \ .
\end{equation}
As $1\leq a\leq2$,
$h_\theta^{(i)}(\BT1\times[0,C])\subseteq\BT1\times[0,1]$, and so we can and
will assume from now on that $C=1$.  Due to the non-continuity of $\varphi^i$
the map $h^{(i)}$ is not continuous, but it still has continuous, strictly
monotonically increasing and strictly concave fiber maps (observe that $C=1$).
For such systems a basic classification was given in \cite{K} (the continuity
assumption made there is not relevant for the facts we are going to state and
use in the following):
\begin{itemize}
\item Ergodic invariant measures correspond to invariant graphs, in
the same sense as in (\ref{eq:associatedmeasures}).
\item There are at most two invariant graphs, one of which is the
0-line. If the Lyapunov exponent of the 0-line is non-positive then
this is the only invariant graph, if its Lyapunov exponent is
positive there exists exactly one other invariant graph $\rho^i$,
which has a negative Lyapunov exponent.
\end{itemize}
Obviously, the Lyapunov exponent of the 0-line in the system $h^{(i)}$
is equal to $\lambda(\mutil^i)$, such that in our situation there always
exists one more $h^{(i)}$-invariant graph $\rho^i$. By
\begin{equation}
\label{eq:gammagraphs}
\Gamma^i(\theta) \ = \ (\varphi^i(\theta),\rho^i(\theta))
\end{equation}
we can then define a $\Lambdatil$-invariant graph, which must be the
support of the measure $\nu$. Again, the Lyapunov exponent
$\lambda_{\textrm{rad}}(\nu)$ is equal to the Lyapunov exponent of
$\rho^i$ for
the
system $h^{(i)}$, and therefore strictly negative.

Summarizing we have found that there exist exactly four ergodic invariant
measures for $\Lambdatil$: $\tilde\mu^s$ and $\tilde\mu^u$, which are embedded
in the 0-torus $S=\BT2 \times \{0\}$ and have positive radial Lyapunov
exponents, and two measures $\nu^s$ and $\nu^u$ which are associated to the
$\Lambdatil$-invariant graphs $\Gamma^s$ and $\Gamma^u$ and have negative
radial Lyapunov exponents.  Among these four measures only $\nu_s$ has
negative exponents in the base, i.e. in $\alpha$-direction, and also in radial
direction.

\ \\
\textit{The global and one-point attractor for $\Lambdatil$.} The
0-torus $S$ is a compact $\Lambdatil$-invariant set, and all ergodic
invariant measures supported on this set have strictly positive
vertical Lyapunov exponents. Therefore it follows from the Uniform
Ergodic Theorem (in fact from a slight generalization, see
\cite{sturman/stark:2000}) that some iterate of $\Lambdatil$ is
uniformly expanding in the vertical direction on a neighborhood of
$S$. Consequently, for sufficiently small $\epsilon$ and suitable
$n\in\BN$ we have $\Lambdatil^n(\BT2 \times [\epsilon,1]) \subseteq
\BT2 \times (\epsilon,1]$. Let $K := \BT2 \times [\epsilon,1]$.  $K$
is compact and forward invariant, and all ergodic invariant measures
supported on $K$ (namely $\nu^s$ and $\nu^u$) have strictly negative
Lyapunov exponents.  Therefore the convergence of the ergodic limits
is again uniform, and a suitable iterate of $\Lambdatil$ is a uniform
vertical contraction on $K$.  But this implies immediately that ${\cal
T} := \bigcap_{n\in\BN} \Lambdatil^n(K)$ is homeomorphic to the
driving space $\BT2$, i.e.\ can be represented as the graph of a
continuous function $T : \BT2 \rightarrow [\epsilon,1]$ (in fact $T$
will be H\"older continuous, see \cite{stark:1999}). Evidently
${\cal T}$ is the boundary of the global attractor ${\cal G}$ and for
all $(\Theta,r) \in Y \setminus S$ there holds
\begin{equation}
\label{eq:pointconv}
|\pi_3(\Lambdatil^n(\Theta,r)) - T(f^n\Theta)| \ \rightarrow \ 0 \ \ \ (n
\rightarrow \infty) \ .
\end{equation}
The one-point attractor mentioned in the beginning is the graph
$\Gamma^u$: The fact that it attracts Lebesgue-a.e.\ initial condition
follows from the fact that on the base this is true for the graph
$\varphi^u$, and in the additional third coordinate the convergence is
given by (\ref{eq:pointconv}). Figure~\ref{fig:1} shows the graph of
$\Gamma^u$ embedded in the manifold ${\cal T}$ and also the graph of
$\varphi^u$, its projection to the 2-dimensional base.

Finally note that $\Lambdatil$ has \sdic, both on the whole phase
space and on ${\cal A} = \textrm{cl}(\Gamma^u)$. For the whole phase
space, this follows from the fact that the base map already has \sdic\
on $\BT2$ by Lemma~\ref{lem:nonhomotopic}~. On the other hand, as the
attractor
is
embedded in ${\cal T}$, the dynamics of $\Lambdatil_{|{\cal A}}$ are
equivalent to the dynamics of $f_{|\textrm{cl}(\Phi^u)}$ (note that
$\textrm{cl}(\Phi^u) = \pi({\cal A})$). Therefore \sdic\ on ${\cal
   A}$ follows either from Proposition~\ref{prop:circle} (if
$\textrm{cl}(\Gamma^u) \neq \BT2$) or
Proposition~\ref{lem:nonhomotopic} (if ${\cal A} = \BT2$).

\ \\
\textit{The original system $\Lambda$.}  Now we can use the results on
$\Lambdatil$ to describe the dynamics of its extension $\Lambda$. The
preimage of $\Gamma^u$ under the factor map $h$ is invariant, consists
of exactly two points on every fiber and attracts Lebesgue-a.e.\
initial condition. As mentioned, the only question we have to leave
open here is whether this two-point attractor further (measurably)
decomposes into two one-valued invariant graphs.

As $\Lambdatil$, the map $\Lambda$ has \sdic\ on the whole phase space
and on the attractor. For the attractor this is immediate as it is
embedded in the two-dimensional torus ${\cal T}_0 := h^{-1}({\cal T})$
such that the dynamics on ${\cal T}_0$ are a two-to-one extension of
$f$, and $f_{|\pi({\cal A})}$ has \sdic. On the whole phase space the
only problem is that in a neighborhood of the 0-line the metric on the
factor space $Y'$ is not equivalent to the usual euclidean metric on
$X'$. However, the two metrics are equivalent if we restrict to the
compact and $\Lambda$-invariant set $h^{-1}(K)$, and as any open set
$U$ which is bounded away from the 0-line ends up in $h^{-1}(K)$ after
a finite number of iterates we obtain that $X' \subseteq \SD_\epsilon$
for a suitable $\epsilon >0$. $X \subseteq \SD_\epsilon$ then follows
again from the fact that $SD_\epsilon$ is closed.

\end{document}